\documentclass[final]{proceedings}

\usepackage{epsfig}
\usepackage{graphicx}
\usepackage{makeidx}
\usepackage{multicol}
\usepackage{subeqn}
\usepackage{latexsym}
\usepackage{url}

\usepackage{bbold}

\usepackage{crop}
\crop
\makeindex

\newtheorem{observation}[theorem]{Observation}

\newcommand{\R}{\mbox{$\mathbb{R}$}}
\newcommand{\Q}{\mbox{$\mathbb{Q}$}}

\newcommand{\conv}{\mathop{\rm conv}}

\newcommand{\Qplus}{\Q^{\geq 0}}
\newcommand{\net}[1]{\mathcal{N}\!\left({#1}\right)}
\newcommand{\setdef}[2]{\left\{{#1}\ :\ {#2}\right\}}
\newcommand{\expans}[1]{\mathcal{X}\!\left({#1}\right)}
\newcommand{\NP}{\mbox{\rm NP}}
\newcommand{\verts}[1]{\mbox{\rm vert}\!\left({#1}\right)}
\newcommand{\wall}[2]{W_{#1}\!\left({#2}\right)}
\newcommand{\spwall}[1]{W\!\left({#1}\right)}
\newcommand{\wallbar}[2]{\overline{W}_{#1}\!\left({#2}\right)}
\newcommand{\cube}[1]{C_{#1}}
\newcommand{\graph}[1]{\mathcal{G}\!\left({#1}\right)}
\newcommand{\bip}[1]{\mathcal{B}\!\left({#1}\right)}
\newcommand{\trees}[1]{\mathcal{T}\!\left({#1}\right)}
\newcommand{\prob}[1]{\mathop{prob}\!\left[{#1}\right]}
\newcommand{\GF}[1]{\mathop{GF}\!\left({#1}\right)}

\newcommand{\mirror}[2]{{#1}^{({#2})}}

\begin{document}

\chapter[On the Expansion of Graphs of 0/1-Polytopes]%
{On the Expansion of Graphs of 0/1-Polytopes}

\begin{authorline}
Volker Kaibel\thanks{Sekr.~MA~6--2, Institut f\"ur Mathematik, Fakult\"at~II,
  TU~Berlin, Stra\ss e des 17. Juni~136, 10623~Berlin, Germany,
  \texttt{kaibel@math.tu-berlin.de}, supported by the Deutsche
  Forschungsgemeinschaft, FOR~413/1--1 (Zi~475/3--1).}
\end{authorline}

\paragraph{Abstract.}
The edge expansion of a graph is the minimum quotient of the number of
edges in a cut and the size of the smaller one among the two node sets
separated by the cut. Bounding the edge expansion from below is
important for bounding the ``mixing time'' of a random walk on the
graph from above. It has been conjectured by Mihail and Vazirani
(see~\cite{FM92}) that the graph of every 0/1-polytope has edge
expansion at least one. A proof of this (or even a weaker) conjecture
would imply solutions of several long-standing open problems in the
theory of randomized approximate counting.  We present different
techniques for bounding the edge expansion of a 0/1-polytope from
below. By means of these tools we show that several classes of
0/1-polytopes indeed have graphs with edge expansion at least one.
These classes include all 0/1-polytopes of dimension at most five, all
simple 0/1-polytopes, all hypersimplices, all stable set polytopes, and
all (perfect) matching polytopes.

\paragraph{MSC~2000.}
52B12 (special polytopes),
52B11 ($n$-dimensional polytopes), 
52B05 (combinatorial properties), 
68W20 (randomized algorithms),
60G50 (random walks)

\section{Introduction}

In the early days of polyhedral combinatorics there was some hope that
investigations of the graphs of 0/1-polytopes that are associated with
certain sets of combinatorial objects might yield insights that could
be exploited in designing algorithms for related combinatorial
optimization problems.  Certainly this hope was inspired by the
success of Dantzig's simplex algorithm for linear programming.  Quite
soon, people came across astonishing facts like the one that the
diameter of the asymmetric traveling salesman polytope equals one for
at most five cities and two for more than five cities (Padberg and
Rao~\cite{PR74}, apparently already discovered, but not published, in
the early fifties by Kuhn~\cite{Kuh01}). This was even outperformed by
the cut polytope of the complete graph on~$n$ nodes that has diameter
one for all~$n\geq 2$ (Barahona and Mahjoub~\cite{BM86}).  Other
polytopes turned out to have more complicated graphs, e.g., the stable
set polytopes, for which two vertices are adjacent if and only if the
symmetric difference of the corresponding stable sets induces a
connected graph~\cite{Chv75}.  Another interesting example is the
basis polytope of a matroid (i.e., the convex hull of the
characteristic vectors of its bases), where two vertices are adjacent
if and only if the corresponding bases have a symmetric difference of
cardinality two (observed by Edmonds in the early 1970's).  All in all
lots of interesting results on the graphs of \emph{special}
0/1-polytopes have been obtained---however, usually without much
impact on algorithms for related optimization problems.

Maybe the best-known result on graphs of \emph{general} 0/1-polytopes
is due to Naddef. He proved~\cite{Nad89} that the graph of any
$d$-dimensional 0/1-polytope has diameter at most~$d$, and thus,
0/1-polytopes satisfy the Hirsch conjecture (claiming that the graph
of any $d$-dimensional polytope with~$n$ facets has diameter at most
$n-d$). Some results on cycles of the graphs of general 0/1-polytopes
have been proved as well by Naddef and Pulleyblank in the
1980's~\cite{Nad84,NP84}. Nevertheless, the graphs of (general)
0/1-polytopes did not receive too much attention. Probably this was due
to the fact that people did not see how to exploit potential knowledge
on this topic with respect to algorithms for combinatorial
optimization problems, where the interest in 0/1-polytopes originally
came from. As for a source of general results on  0/1-polytopes we
refer to~\cite{Zie00}.

The question on graphs of 0/1-polytopes treated in this paper is
mainly motivated by the goal to design algorithms that generate random
elements in classes of combinatorial objects, which often translates
to the task of generating random vertices of 0/1-polytopes. Of course,
in general this includes combinatorial optimization problems via
appropriate choices of random distributions, but here, we will be more
concerned with the task of drawing a vertex according to the uniform
distribution. Maybe the most important motivation of generating
(uniformly distributed) random elements from a set of combinatorial
objects is the fact that in many cases this allows to count the
number of objects approximately by a randomized algorithm. The first
spectacular success of this method was Jerrum and Sinclair's
randomized approximation algorithm for computing the permanent in a
certain large class of 0/1-matrices~\cite{JS89} (extended to arbitrary
matrices with nonnegative integer entries by Jerrum, Sinclair, and
Vigoda~\cite{JSV00}).

For an introduction into the topic of randomized approximate counting
and random generation see~\cite{JS97} or~\cite{Beh99}. Here we briefly
sketch the ideas on the example of the spanning trees of a given graph,
although the exact number of spanning trees can be computed efficiently
by Kirchhoff's matrix tree theorem (see, e.g., \cite[Chap.~24]{AZ01}).

Let~$\trees{G}$ be the set of spanning trees of a graph~$G$. The basic
idea for counting spanning trees via generating them randomly is the
following. Suppose, $G'$ is the graph~$G$ plus an additional
edge~$e'$, and assume, that we do already know a number~$\tau'$
approximating~$|\trees{G'}|$. If we generate a large
set~$\mathcal{T}'$ of spanning trees in~$G'$ uniformly at random, and
if~$\alpha$ is the fraction of those trees in~$\mathcal{T}'$ that do
not contain~$e'$, then we might hope that~$|\trees{G}|$ approximately
equals~$\alpha\cdot\tau'$. Since the number of spanning trees of the
complete graph on~$n$ nodes is well-known to be~$n^{n-2}$, this
suggests an iterative method to approximately compute~$|\trees{G}|$ by
a randomized algorithm.

We do not go into the details of this algorithm and its analysis, but
rather turn to the question how to generate a spanning tree in a
graph~$G$ uniformly at random, where our exposition here is just meant
to give an idea of the method as far as it is useful for understanding
the motivation of the questions on 0/1-polytopes we will consider in
this paper.  The strategy is to perform a (finite) random walk on the
set~$\trees{G}$, meaning that one starts with an arbitrary spanning
tree~$T_0\in\trees{G}$, slightly modifies~$T_0$ randomly to a spanning
tree~$T_1$, slightly modifies~$T_1$ randomly to~$T_2$, and so on.
After a certain number of steps one stops and takes the current tree
as the desired random object. The passage from~$T_i$ to~$T_{i+1}$
could be performed in the following way.  For technical reasons, we
first flip an unbiased coin in order to decide if we ``do nothing''
and stay at $T_{i+1}:=T_i$, or if we try to get to a modified tree as
described subsequently. We first choose a pair $(e,f)$ of edges of~$G$
uniformly at random. If it happens that $e\not\in T_i$ and $f$ lies on
the cycle in~$T_i\cup\{e\}$ then we proceed
to~$T_{i+1}:=T\setminus\{f\}\cup\{e\}$. Otherwise, we stay at
$T_{i+1}:=T_i$.

\begin{figure}[ht]
  \begin{center}
    \includegraphics[height=4cm]{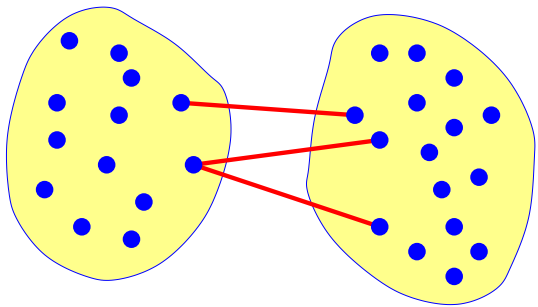}
    \caption{If the neighborhood structure on which a random walk is performed allows to partition the objects into two large parts with only a few connections between them, then the random walk  cannot converge quickly. Fortunately, the converse of this statement is true as well.}
    \label{fig:expans}
  \end{center}
\end{figure}

Thus, we perform a \emph{random walk} in the graph
$\mathcal{G}(\trees{G})$ that has the spanning trees of~$G$ as its
nodes, where two trees are connected if and only if their symmetric
difference consists of two edges. All \emph{transition probabilities}
(i.e., for each ordered pair~$T$ and~$T'$ of adjacent nodes
in~$\mathcal{G}(\trees{G})$ the probability that we proceed to~$T'$ if
we currently are at~$T$) equal $\frac{1}{2\cdot m^2}$, where~$m$ is
the number of edges of~$G$. By standard arguments (see
Section~\ref{sec:expeig}) one can prove that the random walk will be at
each spanning tree with the same probability at step~$i$ if~$i$ tends
to infinity, no matter at which spanning tree we started.  However,
for algorithmic purposes it is of course important that this
convergence does not happen too slow. Responsible for the speed of
convergence is the \emph{edge expansion} of~$\mathcal{G}(\trees{G})$ (see
Figure~\ref{fig:expans}), where the edge expansion of a graph~$H=(V,E)$ is
the number
\begin{eqnarray}
\expans{H}
  &:=& \min\setdef{\frac{|\delta(S)|}{\min\{|S|,|V\setminus
                          S|\}}}{S\subset V,\ S\not=\emptyset,V}\nonumber\\
  & =& \min\setdef{\frac{|\delta(S)|}{|S|}}{S\subset V,\ S\not=\emptyset,|S|\leq\frac{|V|}{2}}\nonumber
\end{eqnarray}
(with~$\delta(S)$ denoting the set of all edges with one end node
in~$S$ and the other one in~$V\setminus S$).
If~$\expans{\mathcal{G}(\trees{G})}$ is bounded by the reciprocal of a
polynomial in the size of~$G$, then the random walk described above
converges ``sufficiently fast.'' Actually, it is well-known that in
our case even~$\expans{\mathcal{G}(\trees{G})}\geq 1$ holds (see the
remarks at the end of Section~\ref{sec:flowmeth}).

Viewing this example of generating spanning trees randomly as a
prototype, one might formulate a strategy for random generation of
certain combinatorial objects as follows. First, one has to choose a
neighborhood structure on the objects and then, transition
probabilities have to be assigned appropriately. Here,
``appropriately'' means (a) that the random walk should asymptotically
behave according to the desired probability distribution and (b) it
should do so approximately already after a small number of steps. Let
us assume that the distribution we aim at is the uniform distribution.
Then, provided that the neighborhood structure is (as in the example)
symmetric and connected, we can achieve goal~(a) always by choosing
the same probability for all proper transitions. In this case,
goal~(b) is equivalent to choosing a neighborhood structure with a
``not too small'' edge expansion.  Of course, in order to be able to
efficiently simulate the random walk it should be also possible to
draw for each object uniformly at random one of its neighboring
objects. However, this will not be at our focus here.

Thus, we are faced with the task to come up with good candidates for
neighborhood structures. Suppose that the set of objects we are
interested in is a family of subsets of a finite set (like in the
example of spanning trees). Then the graph of the associated polytope
(the convex hull of the characteristic vectors of the subsets in the
family) is a natural candidate, where the graph is defined by the
$1$-skeleton, i.e., the zero- and the one-dimensional faces. In fact,
the neighborhood structure we considered in the example is given by
the graph of the spanning tree polytope. Two vertices of that polytope
are adjacent if and only if the symmetric difference of the
corresponding spanning trees consists of two edges (since the spanning
trees of some graph are the bases of a matroid, the \emph{graphic}
matroid defined by that graph).

As mentioned above, two vertices of a stable set polytope are adjacent
if and only if the symmetric difference of the corresponding stable
sets induces a connected subgraph~\cite{Chv75}. Since matchings
correspond to stable sets in the line graph, two vertices of a
matching polytope thus are adjacent if and only if the symmetric
difference of the corresponding matchings is connected. The same is
true for perfect matching polytopes, since they are faces of matching
polytopes. Two of the most prominent random walks in combinatorics are
the ones designed and analyzed by Jerrum and Sinclair~\cite{JS89} on
the set of (near-)perfect matchings of a bipartite graph and on the
set of all matchings of an arbitrary graph.  While the first one lead
to a randomized approximation algorithm for the permanent (for a
certain class of 0/1-matrices), the second one yielded a randomized
approximation algorithm for evaluating the partition function of a
monomer-dimer system in statistical physics, which is the same as the
generating function of the matchings in an arbitrary graph.  In both
cases, the random walk was performed on a subgraph of the graph of the
associated 0/1-polytope, and the crucial step was to prove that this
subgraph has a large edge expansion.  Another example is the random
generation of 0/1-knapsack solutions (leading to a randomized
approximation algorithm for counting as well) due to Morris and
Sinclair~\cite{SM99}. The key step in their result again was to show
that a certain subgraph of the graph of the 0/1-knapsack polytope has
large edge expansion.

It seems to be clear from these examples that it is important to
investigate the question for the edge expansion of general 0/1-polytopes
(i.e., the convex hulls of arbitrary sets of points with coordinates
from~$\{0,1\}$).  Actually, it appears from a citation in a paper of
Feder and Mihail~\cite{FM92} (which we will be concerned with in
Section~\ref{sec:flowmeth}) that Mihail and Vazirani have considered
this question some time ago. Feder and Mihail (and also
Mihail~\cite{Mih92}) quote them with the conjecture that the graph of
every 0/1-polytope has edge expansion at least one. Of course, even a proof
showing that the edge expansion of the graph of any $d$-dimensional
0/1-polytope is bounded by one over a polynomial in~$d$ would be very
important (see also Section~\ref{sec:remarks}).

While this extensive introduction was intended to shed some light on
the relevance of the question for expansion properties of graphs of
0/1-polytopes, the rest of the paper is meant to support the
conjecture of Mihail and Vazirani by some partial results. In
Section~\ref{sec:smalldim} we show that the conjecture indeed is true
for every 0/1-polytope whose dimension does not exceed five. In
Section~\ref{sec:expeig} we list a few well-known facts on random
walks.  The main goal for this is to provide some background that is
relevant for Section~\ref{sec:smalldim}.  As a side effect, the
concepts treated in this introduction may become a bit more clearer.
In Section~\ref{sec:flowmeth} we present some methods for bounding the
edge expansion that are especially suited for graphs of (certain)
0/1-polytopes. In particular, it will turn out that simple
0/1-polytopes, hyper-simplices, and stable set polytopes satisfy
Mihail and Vazirani's conjecture.  We conclude with some remarks in
Section~\ref{sec:remarks}.

The results presented in Sec.~\ref{sec:smalldim} have been obtained in
joined work with Janina Werner~\cite{Wer01}.

\section{Expansion and Eigenvalues}
\label{sec:expeig}

The aim of the present section is to explain the connection between
the edge expansion of a graph and the second largest eigenvalue of a
certain matrix, which will be relevant in Section~\ref{sec:smalldim}. This
connection originates in Alon's and Milman's work~\cite{Alo86,AM86}
and was specifically adapted for our context by Aldous~\cite{Ald87}.
Our treatment closely follows Behrend's book~\cite{Beh99}.

Let~$G=(V,E)$ be a graph (without loops or multiple edges) on~$n:=|V|$
nodes. We define a random walk (i.e., transition probabilities for all
edges---in both directions) on~$G$ in a canonical way.
Let~$\Delta_{\max}$ be the maximum degree of a vertex in~$G$. Each
pair $(v,w)$ of vertices such that $\{v,w\}\in E$ is an edge of~$G$
receives a constant transition probability
$p_{vw}:=\tau:=\frac{1}{2\cdot\Delta_{\max}}$. If~$v\in V$ is a node
of degree~$\Delta_v$, then we set
$p_{vv}:=\frac{1}{2}+(\Delta_{\max}-\Delta_v)\cdot\tau$.
Let~$P\in\R^{V\times V}$ be the matrix with entries $p_{vw}$ ($v,w\in
V$). As defined here, $P$ is a symmetric doubly-stochastic matrix with
a real spectrum
$\lambda_1=1>\lambda_2\geq\lambda_3\geq\dots\geq\lambda_n\geq 0$.  Let
$M\in\R^V$ be a matrix whose columns are eigenvectors of~$P$ that form an
orthonormal basis of~$\R^V$ such that the $i$-th column is an eigenvector
for the eigenvalue~$\lambda_i$. In particular, the first column of~$M$ is
$(\frac{1}{\sqrt{n}},\dots,\frac{1}{\sqrt{n}})$.
Then we have
$$
P=M\cdot\left(
  \begin{array}{cccc}
    1      & 0         & \dots  & 0        \\
    0      & \lambda_2 & \dots  & 0        \\
    \vdots & \vdots    & \ddots & \vdots   \\
    0      & 0         & \dots  & \lambda_n
  \end{array}
  \right)\cdot M^T
$$
(after suitably numbering the vertices of~$G$). 

If the row vector $\pi\in\R^V$ describes the probability distribution
for the start vertex of the random walk, then the distribution after
performing~$i$ steps of our random walk is given by $\pi\cdot P^i$,
i.e., by
\begin{equation}
\label{eq:diag}
\pi\cdot M\cdot\left(
  \begin{array}{cccc}
    1      & 0           & \dots  & 0          \\
    0      & \lambda_2^i & \dots  & 0          \\
    \vdots & \vdots      & \ddots & \vdots     \\
    0      & 0           & \dots  & \lambda_n^i
  \end{array}
  \right)\cdot M^T\enspace.
\end{equation}
For $i\longrightarrow\infty$ this converges to
$$
\pi\cdot M\cdot\left(
  \begin{array}{cccc}
    1      & 0      & \dots  & 0          \\
    0      & 0      & \dots  & 0          \\
    \vdots & \vdots & \ddots & \vdots     \\
    0      & 0      & \dots  & 0
  \end{array}
  \right)\cdot M^T
=
\pi\cdot\left(
  \begin{array}{cccc}
    \frac{1}{n}      & \frac{1}{n}      & \dots  & \frac{1}{n} \\
    \frac{1}{n}      & \frac{1}{n}      & \dots  & \frac{1}{n} \\
    \vdots             & \vdots             & \ddots & \vdots        \\
    \frac{1}{n}      & \frac{1}{n}      & \dots  & \frac{1}{n}
  \end{array}
\right)
=\left(\frac{1}{n},\dots,\frac{1}{n}\right)
\enspace.
$$

Thus, as it was intended, asymptotically the random walk will give
convergence to the uniform distribution over~$V$, independently of the
start distribution (e.g., independent from the start vertex).
Moreover, it follows from~(\ref{eq:diag}) that the speed of
convergence is determined by the second largest
eigenvalue~$\lambda_2$. Intuitively it seems to be clear that the
edge expansion of~$G$ determines how fast the convergence happens. And, in
fact, there is the following strong connection between the edge expansion
and~$\lambda_2$ (see~\cite[Theorem~11.3]{Beh99}).

\begin{theorem}
\label{thm:evexp}
  Let~$G$ be a graph with maximum degree~$\Delta_{\max}$, and
  let~$0\leq\lambda_2<1$ be the second largest eigenvalue of the
  matrix~$P$ defined as above. Then we have
  $$
  (1-\lambda_2)\cdot\Delta_{\max}
  \ \leq\ \expans{G}\ \leq\
  \sqrt{8\cdot(1-\lambda_2)}\cdot\Delta_{\max}\enspace.
  $$

\end{theorem}

The original application of this theorem was, of course, to derive
upper bounds on the size of~$\lambda_2$ by the edge expansion, since the
latter one seems to be easier to access in structural analyses than
the first one. However, with respect to algorithmic issues the
situation is somehow the other way around. While computing the
edge expansion is $\NP$-hard (see Theorem~\ref{thm:nphard}), the second
largest eigenvalue can be calculated efficiently. We will exploit this
fact in the next section.

\section{Small Dimensions}
\label{sec:smalldim}

Aichholzer classified all 0/1-polytopes of dimension less than or
equal to five up to isometries of the cube, i.e., up to flipping and
permuting the coordinates~\cite{Aic00}.  Table~\ref{tab:01polys}
shows the number of classes for each dimension.
\begin{table}[ht]
  \caption{The numbers of classes of 0/1-polytopes.}
  \begin{center}
    \begin{tabular}{l|rrrrr}
      Dimension  & 1 & 2 &  3 &   4 &       5 \\
      \hline
      \# Classes & 1 & 2 & 12 & 349 & 1226525
    \end{tabular}
    \label{tab:01polys}
  \end{center}
\end{table}

Thus, in principle one can compute the edge expansion of the graph of each
0/1-polytope up to dimension five by computer. Unfortunately, the
following result shows that in general, computing the edge expansion is
difficult. This is well-known for some time (e.g., \cite{LR88}).
However, since we could not find an explicit proof in the literature,
we include one here.

\begin{theorem}
  \label{thm:nphard}
  The problem of computing~$\expans{G}$ for arbitrary graphs~$G$ is
  $\NP$-hard.
\end{theorem}

\begin{proof}
  We reduce the problem of finding a maximum (unweighted) cut in a
  graph (which was proved to be $\NP$-hard by Karp~\cite{Kar72}) to
  the problem of computing the edge expansion of some related graph.
  The proof is an extension of the proof of the $\NP$-hardness of the
  equicut problem given by Garey, Johnson, and
  Stockmeyer~\cite{GJS76}.
  
  Let~$G=(V,E)$ be a graph with $n:=|V|$ nodes. We construct a graph
  $G'=(V',E')$, where $V'=V\uplus W$ for some set~$W$, disjoint
  from~$V$, with $|W|=n$, and with~$E'$ containing all possible edges
  except the ones in~$E$. Thus, $G'$ has $n'=2n$ nodes. We denote by
  $\delta_G(S)$ and $\delta_{G'}(S')$ the set of all edges of~$G$
  respectively~$G'$ having precisely one end node in~$S$
  respectively~$S'$ and
  define
  $$
  \eta_{G'}(S'):=\frac{|\delta_{G'}(S')|}{\min\{|S'|,|V'\setminus
    S'|}\enspace.
  $$

  We first show that it suffices to consider node subsets of
  cardinality~$\frac{n'}{2}=n$ in order to compute the edge expansion
  of~$G'$.  Let~$S\subseteq V$ and $T\subseteq W$ be two sets of nodes
  of~$G'$ with $k:=|S|+|T|\leq n$. We have
  $|\delta_{G'}(S\cup T)|=k\cdot(2n-k)-|\delta_G(S)|$ and
  \begin{equation}
    \label{eq:quot1}
    \eta_{G'}(S\cup T)=2n-k-\frac{|\delta_G(S)|}{k}\enspace.
  \end{equation}
  In particular, if $k=n$ then
  \begin{equation}
    \label{eq:quot2}
    \eta_{G'}(S\cup T)=n-\frac{|\delta_G(S)|}{n}\enspace.
  \end{equation}
  holds. 
  
  We claim that the right hand side of~(\ref{eq:quot2}) is less than
  or equal to the right hand side of~(\ref{eq:quot1}) for each~$1\leq
  k\leq n$. Indeed, this claim is equivalent to
  $$
  n-k+\left(\frac{1}{n}-\frac{1}{k}\right)\cdot|\delta_G(S)|\geq 0\enspace,
  $$
  which follows from
  $$
  |\delta_G(S)|\leq |S|\cdot n\leq k\cdot n\enspace.
  $$

  Thus, we have (where the second equation follows from~(\ref{eq:quot2})) 
  \begin{eqnarray}
    \expans{G'} 
    &=& \min\setdef{\eta_{G'}(S\cup T)}{S\subseteq V, W\subseteq W, 
                                                |S|+|T|=n}\nonumber\\ 
    &=& \min\setdef{n-\frac{\delta_G(S)}{n}}{S\subseteq V, W\subseteq W, 
                                                |S|+|T|=n}\nonumber\\                                               
    &=& n-\frac{\max\setdef{|\delta_G(S)|}{S\subseteq V}}{n}\nonumber\enspace.
  \end{eqnarray}
%  and hence,
%  $$
%  \max\setdef{|\delta_G(S)|}{S\subseteq V}=n\cdot(n-\expans{G'})\enspace.
%  $$
\end{proof}

In view of Theorem~\ref{thm:nphard} we decided first to calculate the
lower bounds on the edge expansion provided by Theorem~\ref{thm:evexp} for
each 0/1-polytope of dimension four and five. And, somewhat
surprising, it turned out that for none of the polytopes this bound was
less than one. Thus, the conjecture of Mihail and Vazirani is true for
0/1-polytopes up to dimension five.

\begin{theorem}
  The graph of each 0/1-polytope of dimension less than or equal to
  five has edge expansion at least one.
\end{theorem}

\begin{figure}[ht]
  \begin{center}
    \includegraphics[height=7cm,angle=270]{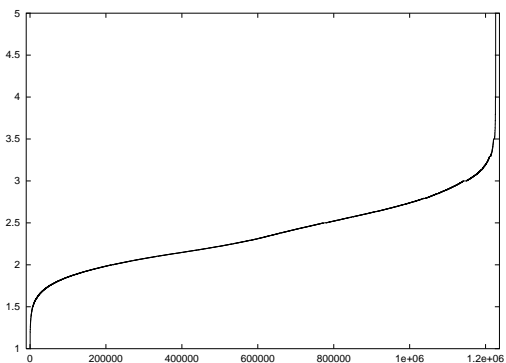}
    \label{fig:histo}
  \end{center}
  \caption{The (lower) eigenvalue bounds on the edge expansion for all 1226525~
    five-dimensional 0/1-polytopes.}
\end{figure}

Figure~\ref{fig:histo} shows that in many cases the lower bound given
by the second largest eigenvalue even was significantly larger than
one.

\section{Flow Methods}
\label{sec:flowmeth}

In this section, we describe methods for proving that a graph has good
edge expansion properties that are specifically suited for graphs of
0/1-polytopes. Applying these methods we will show that the conjecture
of Mihail and Vazirani is true for well-known classes of
0/1-polytopes (see Corollaries~\ref{cor:ssh} and~\ref{cor:cswalls}).
On the other hand it will be quite obvious that the methods are not
sufficient to prove the conjecture in its whole generality.

\subsection{Expansion and flows}
\label{subsec:expflow}

In order to bound the edge expansion of a graph~$G=(V,E)$ from below
we will construct certain flows in the (uncapacitated)
network~$\net{G}=(V,A)$, where~$A$ contains for each edge~$\{u,v\}\in
E$ both arcs~$(u,v)$ and~$(v,u)$.  This strategy dates back to the
method of ``canonical paths'' developed by Sinclair
(see~\cite{Sin93}). The extension to flows was explicitly exploited by
Morris and Sinclair~\cite{SM99}. Feder and Mihail~\cite{FM92} use
random canonical paths, which can equivalently be formulated in terms
of flows.

The crucial idea is to construct for each ordered pair $(s,t)\in
V\times V$ a flow $\phi_{(s,t)}:A\longrightarrow\Qplus$ in the
network~$\net{G}$ sending one unit of some commodity from~$s$ to~$t$.
Let $\phi:=\sum_{(s,t)\in V\times V}\phi_{(s,t)}$ be the sum of all
these flows. By
$$
\phi_{\max}:=\max\setdef{\phi(a)}{a\in A}
$$
we denote the maximal amount of $\phi$-flow on any arc. By
construction of~$\phi$, the total amount $\phi(S:V\setminus S)$ of
$\phi$-flow leaving~$S$ is at least $|S|\cdot(n-|S|)$, where $n=|V|$.
On the other hand, we have $\phi(S:V\setminus S)\leq
\phi_{\max}\cdot|\delta(S)|$. This implies
$|S|\cdot(n-|S|)\leq\phi_{\max}\cdot|\delta(S)|$, and hence, if
$|S|\leq\frac{n}{2}$ holds,
$$
\frac{|\delta(S)|}{|S|}\geq \frac{n}{2\cdot\phi_{\max}}\enspace.
$$
Thus, we have proved
\begin{equation}
  \label{eq:expphimax}
  \expans{G}\geq\frac{n}{2\cdot\phi_{\max}}\enspace.
\end{equation}

In the light of inequality~(\ref{eq:expphimax}) it is clear that the
task is to construct a flow~$\phi$ as above with $\phi_{\max}$ as
small as possible in order to prove a strong lower bound on the
edge expansion of~$G$. 

\subsection{Fractional wall-matchings}
\label{subsec:wallmatch}

While the setting presented so far applies to general graphs, we now
derive a method to construct~$\phi$ in the special situation where~$G$
is the graph of a 0/1-polytope. The method generalizes ideas for
analyzing random walks on the bases-exchange graph of matroids due to
Feder and Mihail~\cite{FM92}.

Let~$P\subset\R^d$ be a 0/1-polytope. A \emph{wall} of~$P$ is the
intersection of~$P$ with any face of the cube
$\cube{d}:=\setdef{x\in\R^d}{0\leq x_i\leq 1 \mbox{ for all }i}\supseteq P$.
Thus, the walls of~$P$ are special faces of~$P$. Usually, we will
identify a wall of~$P$ with its vertices. The faces~$F$ of~$\cube{d}$
are in one-to-one correspondence with the vectors
$\sigma(F)\in\{0,1,\star\}^d$ (and vice versa) via
$$
F=\setdef{x\in\cube{d}}{x_i=\sigma(F)_i
                        \mbox{ for all $i$ with } \sigma(F)_i\not=\star}\enspace.
$$
For a face $F\not=\cube{d}$ of~$\cube{d}$ let
$\mu(F):=\min\setdef{i}{\sigma(F)_i=\star}$ be the ``smallest
direction'' of~$F$. Let~$W$ be a wall of~$P$ and let~$F$ be the
inclusion minimal face of~$\cube{d}$ with $W=P\cap F$. The vector
$\sigma(W):=\sigma(F)$ indicates the components in which all vertices
of~$W$ agree, and~$\mu(W):=\mu(F)$ is the smallest coordinate
direction of any edge of~$W$. We define
$W_0:=\setdef{w\in W}{w_{\mu(W)}=0}$ and $W_1:=\setdef{w\in
  W}{w_{\mu(W)}=1}$, and denote by $\bip{W}$ the bipartite subgraph
of~$\graph{P}$ induced by the two disjoint subsets~$W_0$ and~$W_1$ of
nodes of~$\graph{P}$.

A wall~$W$ of~$P$ is called \emph{initial} if there is
some~$i\in\{0,1,\dots,d\}$ such that $\sigma(W)_j\in\{0,1\}$ for
$1\leq j\leq i$ and $\sigma(W)_j=\star$ for all $j>i$. The following fact
follows immediately from the definitions.

\begin{lemma}
  \label{lem:uniqueinitwall}
  For every edge~$e$ of a 0/1-polytope~$P$ there is a unique initial
  wall~$W$ of~$P$ such that~$e$ is an edge of~$\bip{W}$.
\end{lemma}

Thus the bipartite graphs associated with the initial walls of~$P$
induce a partition of the edges of~$P$.

\begin{figure}[ht]
  \begin{center}
    \includegraphics[height=4cm]{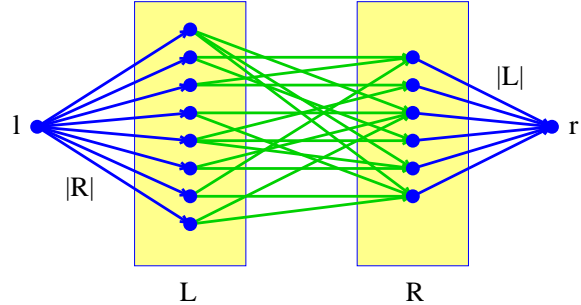}    
    \caption{The bipartite graph on~$L\uplus R$  has a fractional matching if and only
      if in the network indicated in the figure there is a (non
      negative) flow sending~$|L|\cdot|R|$ units of some commodity
      from~$l$ to~$r$.  The arcs leaving~$l$ have capacities~$|R|$,
      the arcs entering~$r$ have capacities~$|L|$, and the arcs
      connecting~$L$ to~$R$ have infinite capacities.}
    \label{fig:fracmatch}
  \end{center}
\end{figure}

A bipartite graph with bipartition $L\uplus R$ has a \emph{fractional
  matching} if one can assign nonnegative weights to its edges such
that all nodes in~$L$ have the same weighted degree, and the same does 
hold for all nodes in~$R$ as well (see Figure~\ref{fig:fracmatch}).

\begin{observation}
\label{obs:fracmatch}
If a bipartite graph~$B$ with bipartition $L\uplus R$ has a fractional
matching and there is a constant amount of some commodity located in
each node in~$L$ (or~$R$, respectively), then one can distribute the
entire amount of the commodity from~$L$ to~$R$ through the edges
of~$B$ such that each node in~$R$ (or~$L$, respectively) receives the
same amount of the commodity.
\end{observation}

\begin{figure}[ht]
  \begin{center}
    \includegraphics[height=3.5cm]{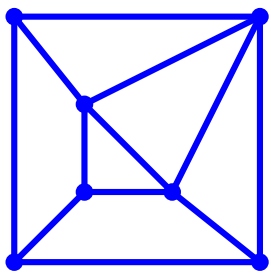}
    \caption{The graph of the 0/1-polytope~$P$ arising from~$\cube{3}$ by
      removing one vertex. Independently of the numbering of the
      coordinate directions $\bip{P}$ has no fractional matching.}
    \label{fig:nofracmatch}
  \end{center}
\end{figure}

A 0/1-polytope~$P$ \emph{has fractional wall-matchings} if $\bip{W}$
has a fractional matching for every wall~$W$ of~$P$.  In general, the
bipartite graph $\bip{W}$ associated to a wall~$W$ of a
0/1-polytope~$P$ does not necessarily have a fractional matching (see
Figure~\ref{fig:nofracmatch}).
However, several interesting classes of 0/1-polytopes have fractional
wall-matchings, as we will show below.  The method to construct
suitable flows~$\phi$ we will describe does only work for such
0/1-polytopes. Thus, from now on we assume that~$P\subset\R^d$ is a
0/1-polytope that has fractional wall-matchings.

\begin{figure}[ht]
  \begin{center}
    \includegraphics[height=3cm]{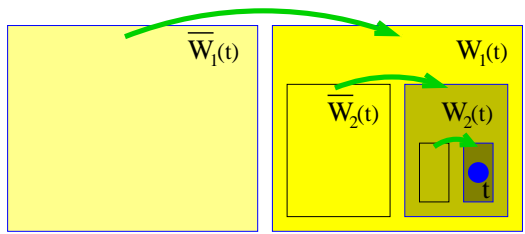}
    \caption{Simultaneous construction of the flows~$\phi_{(s,t)}$ for all $s\in V$.}
    \label{fig:wallflow}
  \end{center}
\end{figure}

Let $t\in\verts{P}$ be a vertex of~$P$. We will particularly be concerned
with the initial walls
$$
\wall{i}{t}:=\setdef{w\in\verts{P}}{w_1=t_1,\dots,w_i=t_i}\qquad(i=0,1,\dots,d)\enspace.
$$
These walls form a \emph{flag} of~$P$, i.e., we have
$$
\{t\}=\wall{d}{t}\subseteq\wall{d-1}{t}\subseteq\dots\subseteq\wall{1}{t}\subseteq\wall{0}{t}=P\enspace.
$$
For each $i\in\{1,\dots,d\}$ we define
$\wallbar{i}{t}:=\wall{i-1}{t}\setminus\wall{i}{t}$. Now we are ready
to construct all flows~$\phi_{(s,t)}$, $s\in V$, simultaneously in~$d$
steps. Imagine a single unit of some commodity initially placed at
each node. Suppose that before we perform step~$i\in\{1,\dots,d\}$
the~$n$ units of the commodity are distributed uniformly among the
nodes in $\wall{i-1}{t}$ (as it is the case before the first step).
Since we have assumed that~$P$ has fractional wall-matchings we can
route (see Observation~\ref{obs:fracmatch}) the amount of commodity
distributed at the nodes in~$\wallbar{i}{t}$ through the arcs
corresponding to the edges of~$\bip{\wall{i-1}{t}}$ such that
afterwards the~$n$ units of our commodity are uniformly distributed
among the nodes in~$\wall{i}{t}$.  Figure~\ref{fig:wallflow}
illustrates the construction.

\begin{figure}[ht]
  \begin{center}
    \includegraphics[height=3.5cm]{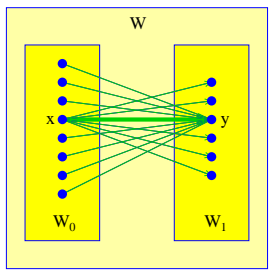}
    \caption{The arc sets $A^-(x)$ and $A^+(y)$.}
    \label{fig:analysis}
  \end{center}
\end{figure}

For each pair $(s,t)\in V\times V$ we thus have defined a
flow  in the network $\net{\graph{P}}$ sending one unit of some
commodity from~$s$ to~$t$.  It remains to
bound the maximal flow~$\phi_{\max}$ produced by $\phi:=\sum_{(s,t)\in
  V\times V}\phi_{(s,t)}$ at any arc. Therefore, let~$(x,y)$ be any arc
of $\net{\graph{P}}$. By Lemma~\ref{lem:uniqueinitwall} there is a
unique initial wall~$W$ of~$P$ such that $\bip{W}$ contains the
edge~$\{x,y\}$ (see Figure~\ref{fig:analysis}).
Due to symmetry reasons we might assume $x\in W_0$ and $y\in W_1$.
Let $A^-(x)$ and $A^+(y)$ be the sets of out-arcs respectively in-arcs
incident to~$x$ respectively~$y$ corresponding to edges of~$\bip{W}$.
In particular, we have $(x,y)\in A^-(x)$ and $(x,y)\in A^+(y)$.  The
arcs going from~$W_0$ to~$W_1$ are only used by the flows
$\phi_{(s,t)}$ with $s\not\in W_1$ and $t\in W_1$. Thus, the total
amount of flow carried by these arcs is
$\frac{|W_0|\cdot|W_1|}{|W|}\cdot n$. Consequently, precisely
$\frac{|W_1|}{|W|}\cdot n$ units of flow are sent through $A^-(x)$ and
$\frac{|W_0|}{|W|}\cdot n$ units of flow are sent through $A^+(y)$.
Hence, $(x,y)$ carries at most
$$
\min\left\{\frac{|W_1|}{|W|},\frac{|W_0|}{|W|}\right\}\cdot n\leq\frac{n}{2}
$$
units of flow. Since this holds for every arc of~$\net{\graph{P}}$,
we have $\phi_{\max}\leq\frac{n}{2}$. By~(\ref{eq:expphimax}) this
proves the following result.

\begin{theorem}
  \label{thm:fracwallmatch}
  If~$P$ is a 0/1-polytope that has fractional wall-matchings, then 
  $\expans{\graph{P}}\geq 1$ holds.
\end{theorem}

Thus, 0/1-polytopes that have fractional wall-matchings satisfy the
conjecture of Mihail and Vazirani.

\subsection{Walls with regular graphs}

\begin{figure}[ht]
  \begin{center}
    \includegraphics[height=3.5cm]{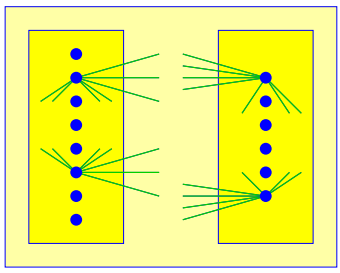}
    \caption{Regular walls yield fractional wall-matchings, since in
      each of the relevant bipartite graphs all vertices in the left shore
      have the same degree and the same is true for all vertices in
      the right shore.}
    \label{fig:regwalls}
  \end{center}
\end{figure}

Let us say that a 0/1-polytope~$P$ has \emph{regular walls} if the
graph of every wall of~$P$ is regular, i.e., all its vertices have the
same degree. It is obvious that every 0/1-polytope
with regular walls has fractional wall-matchings (see
Figure~\ref{fig:regwalls}). This proves the following consequence of
Theorem~\ref{thm:fracwallmatch}. 
\begin{corollary}
\label{cor:regwalls}
  If a 0/1-polytope~$P$ has regular walls then $\expans{\graph{P}}\geq 1$ 
  holds.
\end{corollary}

A $d$-dimensional polytope~$P$ is \emph{simple} if every vertex lies
in precisely~$d$ facets, or, equivalently, if $\graph{P}$ is
$d$-regular.  The polytopes
$$
\conv\setdef{v\in\{0,1\}^d}{\sum_{i=1}^dv_i=\varrho}\qquad(\varrho\in\{0,1,\dots,d\})
$$
are called \emph{hyper-simplices} (they are special Knapsack
polytopes).

\begin{corollary}
\label{cor:ssh}
If a 0/1-polytope~$P$ is simple or a hyper-simplex, then
$\expans{\graph{P}}\geq 1$ holds.
\end{corollary}

\begin{proof}
  Every face of a simple polytope is simple, and thus has a regular
  graph. Every wall of a hyper-simplex is a hyper-simplex, again.
  Since hyper-simplices obviously have a transitive automorphism
  group, they have regular graphs. Thus, in any of the two cases of
  the claim, $\expans{\graph{P}}\geq 1$ holds by
  Corollary~\ref{cor:regwalls}.
\end{proof}

\subsection{Balanced uniform 0/1-polytopes}
\label{subsec:balance}

A 0/1-polytope~$P\subset\R^d$ is called \emph{$\varrho$-uniform}
($\varrho\in\{0,1,\dots,d\}$) if it is contained in the hyperplane
$\setdef{x\in\R^d}{\sum_{i=1}^dx_i=\varrho}$, i.e., if all vertices
of~$P$ have precisely~$\varrho$ ones. For instance, hyper-simplices and
basis polytopes of matroids are uniform.  Obviously, every wall of a
uniform 0/1-polytope is uniform as well.

A 0/1-polytope~$P\subset\R^d$ is \emph{balanced} if for every
$\sigma\in\{0,1,\star\}$ and for each pair $i,j\in\{1,\dots,d\}$ with
$i\not= j$ and $\sigma_i=\sigma_j=\star$ the relation
\begin{equation}
  \label{eq:balance}
  |W_{0,0}|\cdot|W_{1,1}|\leq|W_{1,0}|\cdot|W_{0,1}|  
\end{equation}
holds, where $W$ is the wall of~$P$ defined by~$\sigma$ and
$W_{\alpha,\beta}:=\setdef{w\in W}{w_i=\alpha,w_j=\beta}$.  If
$W_{0,1}\cup W_{1,1}\not=\emptyset$ (i.e., there is some $w\in W$ with
$w_j=1$), then~(\ref{eq:balance}) is equivalent to
\begin{equation}
\label{eq:negcor1}
\frac{|W_{1,1}|}{|W_{0,1}|+|W_{1,1}|}
\leq
\frac{|W_{1,1}|+|W_{1,0}|}{|W_{0,1}|+|W_{1,1}|+|W_{0,0}|+|W_{1,0}|}\enspace.
\end{equation}
This means, that for a vertex~$w$ chosen uniformly at random from~$W$
the probability of the event $w_i=1$ does not increase by conditioning
on the event $w_j=1$. Similarly, (\ref{eq:balance}) is equivalent to
the fact that for a vertex~$w$ chosen uniformly at random from~$W$ the
probability of the event $w_i=0$ does not increase by conditioning on
the event $w_j=0$.

The property of being balanced is not invariant under arbitrary
symmetries of the cube. However, it is invariant under simultaneous
``flipping'' of all coordinates (and under arbitrary permutations of the
coordinates). 

\begin{proposition}
  \label{prop:balanced01}
  Balanced uniform 0/1-polytopes have fractional wall-matchings.
\end{proposition}

We omit the proof, which closely follows the corresponding proof on
the bases-exchange graph of balanced matroids due to Feder and
Mihail~\cite{FM92}.

Proposition~\ref{prop:balanced01} and Theorem~\ref{thm:fracwallmatch}
imply the following.

\begin{theorem}
\label{thm:balance}
  Every balanced uniform 0/1-polytope~$P$ satisfies
  $\expans{\graph{P}}\geq 1$. 
\end{theorem}

A matroid~$\mathcal{M}$ on the ground set~$E$ has the \emph{negative
  correlation property} if for a basis~$B$ chosen uniformly at random
from the set of bases of~$\mathcal{M}$ and for every pair of elements
$e,f\in E$
$$
\prob{e\in B}\geq\prob{e\in B\ |\ f\in B}
$$
holds, i.e., the probability of the event $e\in B$ does not
increase by conditioning on the event $f\in B$. A matroid~$\mathcal{M}$
is \emph{balanced} if every minor of~$\mathcal{M}$ has the negative
correlation property. Regular (in particular: graphic) matroids are
known to be balanced. 

It is obvious that the basis polytope
$P(\mathcal{M}):=\conv\setdef{\chi(B)}{B\mbox{ basis of }\mathcal{M}}$
(where $\chi(B)$ is the characteristic vector of $B\subseteq E$) of a
balanced matroid~$\mathcal{M}$ is uniform and balanced. Thus,
Theorem~\ref{thm:balance} immediately yields
$\expans{\graph{P(\mathcal{M})}}\geq 1$.  Notice that the actual
adjacency structure on $P(\mathcal{M})$ is irrelevant for this.

Hence, Theorem~\ref{thm:balance} generalizes the result of Feder and
Mihail~\cite{FM92} saying that the bases-exchange graphs (where two
bases are adjacent if and only if their symmetric difference has two
elements) of balanced matroids have edge expansion at least one.

\subsection{Cube-spanned walls}
\label{subsec:cswalls}

The technique described in this subsection is particularly suited for
proving that 0/1-polytopes coming from certain combinatorial problems
have graphs with large edge expansion (see Cor.~\ref{cor:cswalls}). It
relies on the high symmetry of the graph~$\graph{Q}$ of a cube~$Q$,
from which one easily derives the following fact (where the
\emph{antipodal vertex} of some vertex~$x$ of~$Q$ is the vertex with
maximum distance from~$x$ in~$\graph{Q}$)

\begin{observation}
  \label{obs:cswalls}
  For a cube~$Q$ it is possible to define for each pair $(s,t)$ of
  \emph{antipodal} vertices a flow~$\psi_{(s,t)}$
  in~$\net{\graph{Q}}$ sending one unit of some commodity from~$s$
  to~$t$ such that for the total flow $\psi:=\sum_{s,t}\psi_{(s,t)}$
  one has $\psi(a)=1$ for each arc~$a$ in~$\net{\graph{Q}}$.
\end{observation}

Let~$P\subset\R^d$ be any 0/1-polytope.  A subset
$C\subseteq\verts{P}$ of vertices of~$P$ is called an \emph{affine
  cube in~$P$} if~$C$ is affinely isomorphic to $\{0,1\}^k$ for
some~$k$, or, equivalently, if there is a
subset~$I\subseteq\{1,\dots,d\}$ (with $|I|=k$) such that the
orthogonal projection of~$\R^d$ onto~$\R^I$ induces a bijection
between~$C$ and $\{0,1\}^I$.  It is not hard to see
that~$C\subseteq\verts{P}$ is an affine cube (in~$P$) if and only if
there are 0/1-vectors $z^{(1)},\dots,z^{(k)}\in\{0,1\}^d$, pairwise
orthogonal to each other, such that for every $x\in C$
$$
C=\setdef{x\oplus\epsilon_1z^{(1)}\oplus\cdots\oplus\epsilon_kz^{(k)}}{\epsilon\in\{0,1\}^k}
$$
(where $\oplus$ denotes addition modulo two). The vertex $x\oplus
z^{(1)}\oplus\cdots\oplus z^{(k)}$ is the \emph{antipodal vertex
  of~$x$ in~$C$}.  In particular, we will be interested in
\emph{affine edge-cubes} in~$P$, i.e., affine cubes in~$P$ on
which~$\graph{P}$ induces the graph of a cube.

\begin{figure}[ht]
  \begin{center}
    \includegraphics[height=3.5cm]{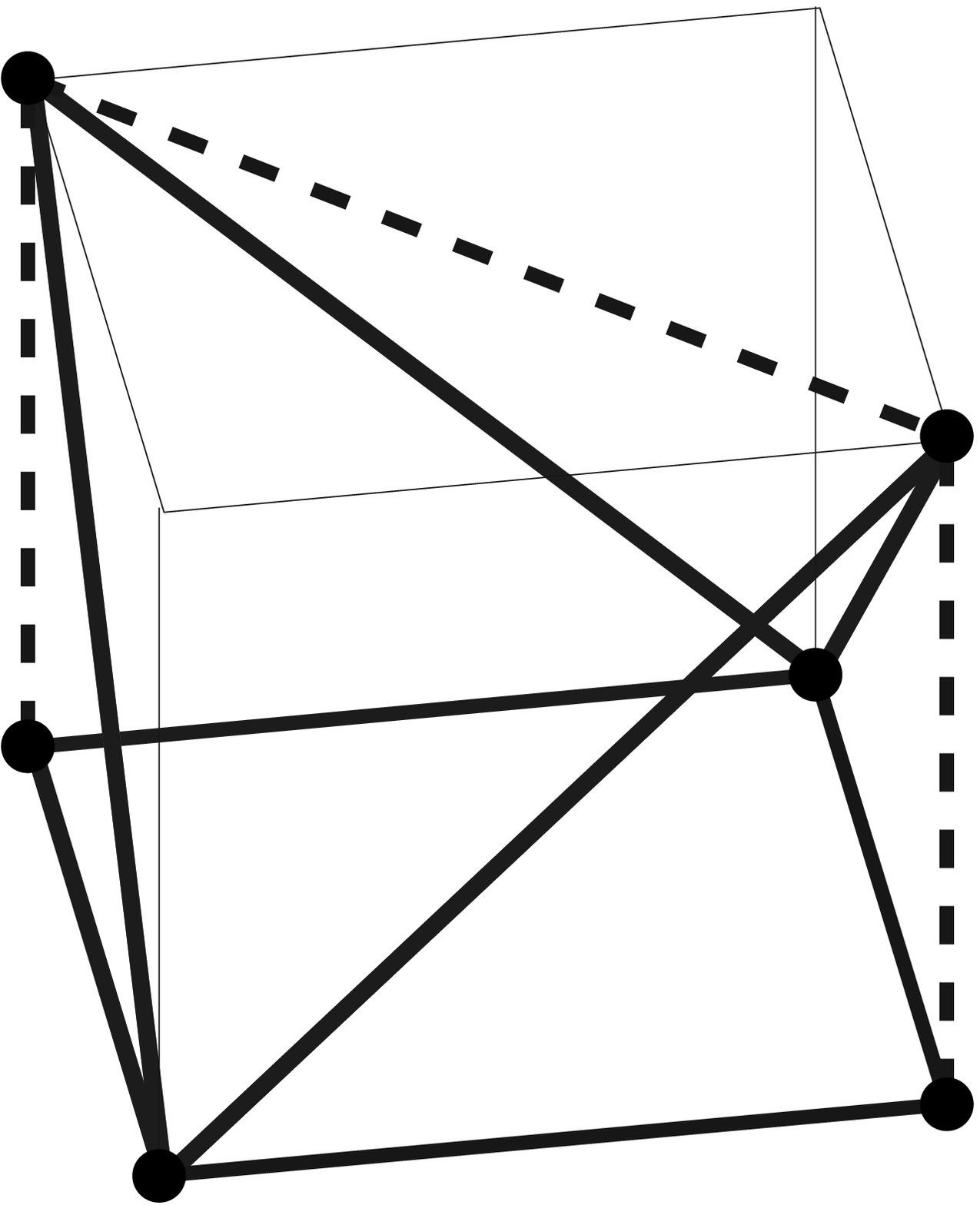}\hfill
    \includegraphics[height=3.5cm]{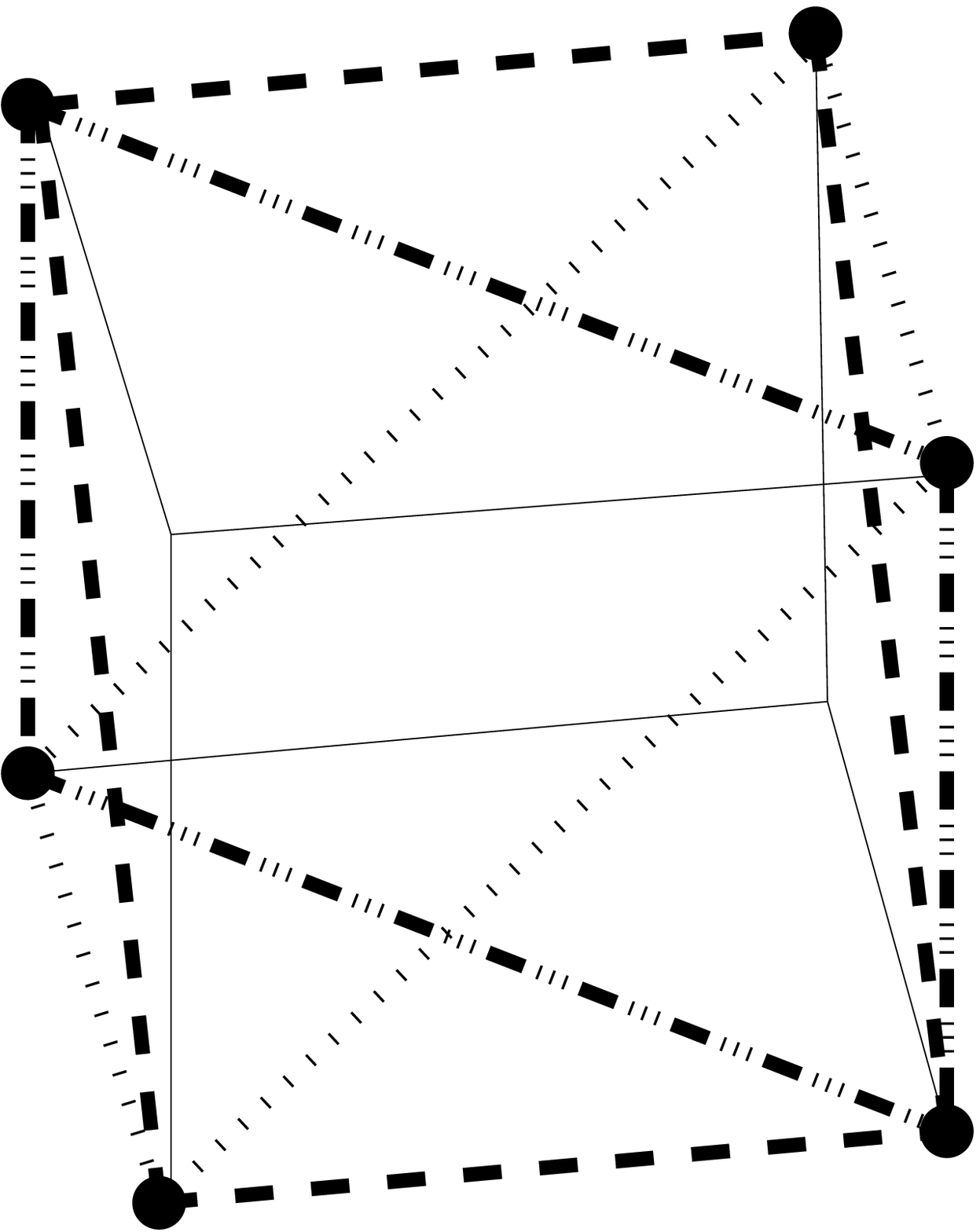}\hfill
    \includegraphics[height=3.5cm]{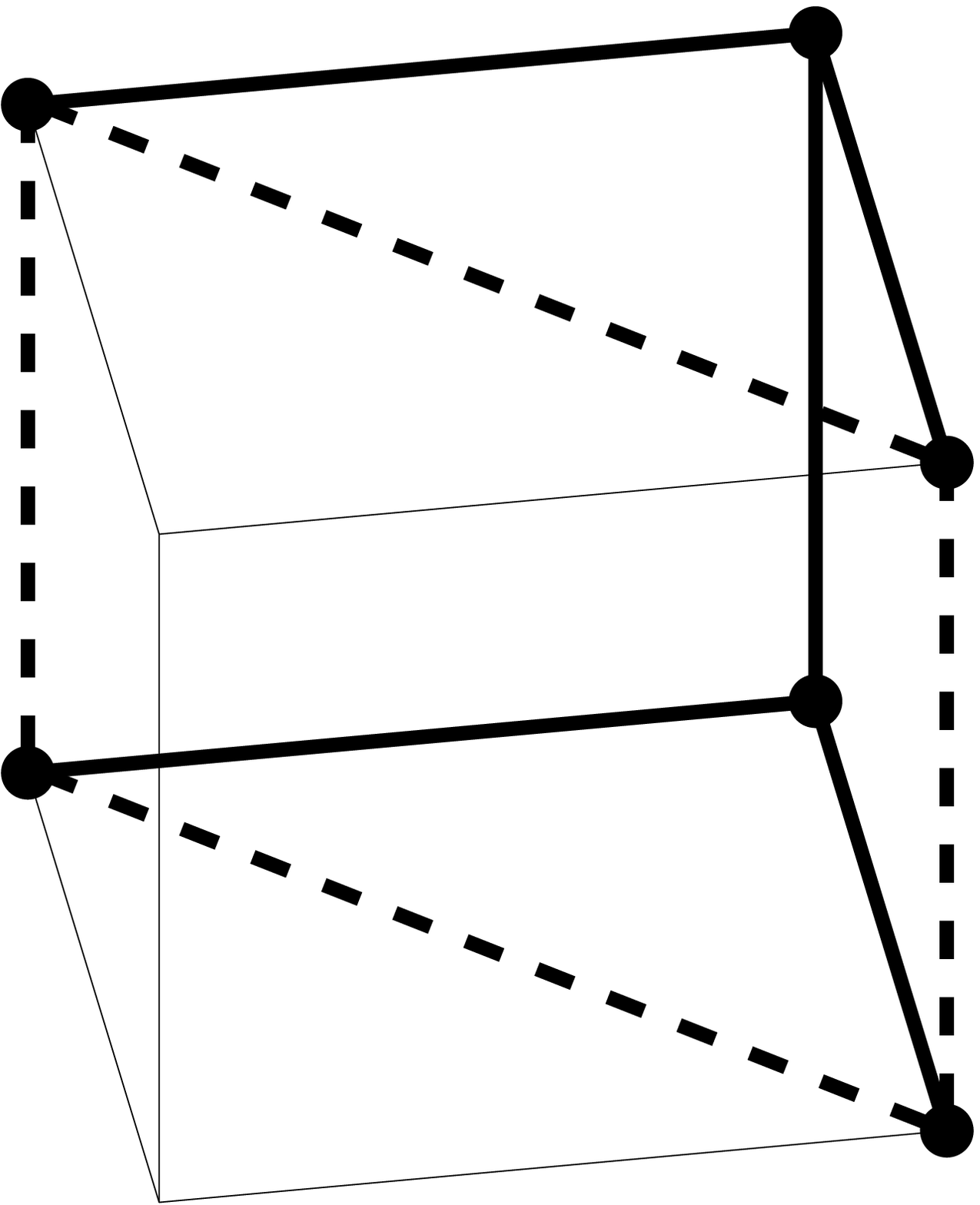}
    \caption{Three $3$-dimensional walls. The first one is spanned by
      a cube (but not edge-cube spanned), the second one is edge-cube
      spanned (but not uniquely edge-cube spanned), and the third one
      is uniquely edge-cube spanned.}
    \label{fig:cs}
  \end{center}
\end{figure}

For a subset $A\subseteq\verts{P}$ let us call the intersection of all
walls of~$P$ that contain~$A$ the \emph{wall spanned by~$A$}, denoted
by~$\spwall{A}$. A wall~$W$ of~$P$ is \emph{edge-cube spanned} if
there is an affine edge-cube~$C$ with $\spwall{C}=W$ (which is
equivalent to the fact that each pair of antipodal vertices in~$C$
spans~$W$). A wall~$W$ of~$P$ is \emph{uniquely edge-cube spanned} if
it is spanned by an affine edge-cube~$C$ in~$P$ and if it is not
spanned by any other affine edge-cube~$C'\not= C$ in~$P$. In this
case, we call the vertices in~$C$ the \emph{cube vertices of~$W$}. See
Figure~\ref{fig:cs} for examples.

For a vertex~$x\in W$ in a wall~$W$ of~$P$ the vertex
$\mirror{x}{W}:=x\oplus t^{(W)}$ is the \emph{mirror image} of~$x$ with
respect to~$W$, where $t^{(W)}$ is the 0/1-vector having ones
precisely in those components where~$\sigma(W)$ has stars. In general,
$\mirror{x}{W}$ needs not to be contained in~$W$. If, however, a wall~$W$ is
spanned by an affine cube~$C$, then ~$\mirror{x}{W}$ is the antipodal vertex
of~$x$ in~$C$ for every $x\in C$; in particular, $\mirror{x}{W}\in W$.

\begin{lemma}
  \label{lem:numwalls} 
  Let~$P$ be a 0/1-polytope, and let $u,v\in\verts{P}$, $u\not= v$, be
  two distinct vertices of~$P$. There are at most
  $\frac{1}{2}|\verts{P}|$ walls~$W$ of~$P$ such that~$u$, $v$, and
  their mirror images with respect to~$W$ are contained in~$W$.
\end{lemma}

\begin{proof}
  Let~$\mathcal{W}(u,v)$ be the set of all walls~$W$ of~$P$ such
  that~$u$, $v$, and their mirror images~$\mirror{u}{W}$,
  $\mirror{v}{W}$ are contained in~$W$.  We have
  \begin{equation}
    \label{eq:numwalls}
    u\oplus\mirror{v}{W}, v\oplus\mirror{u}{W}\leq t^{(W)}
    \qquad\mbox{and}\qquad
    u\oplus\mirror{u}{W}=v\oplus\mirror{v}{W}=t^{(W)}
  \end{equation}
  (where $\leq$ is meant to hold component-wise).  Since $u\not= v$ we
  have $\mirror{u}{W}\not=\mirror{v}{W}$. Thus, we can define a map
  $\omega$ assigning to each $W\in\mathcal{W}(u,v)$ the two-element
  subset $\omega(W):=\{\mirror{u}{W},\mirror{v}{W}\}$
  of~$\verts{P}$.
  
  Suppose, for $W,W'\in\mathcal{W}(u,v)$ we have
  $x\in\omega(W)\cap\omega(W')$. After possibly interchanging the
  roles of~$u$ and~$v$, by~(\ref{eq:numwalls}) we have $v\oplus x\leq
  u\oplus x$ and thus $t^{(W)}=u\oplus x=t^{(W')}$, yielding~$W=W'$.
  Thus the images of~$\omega$ have pairwise empty intersections, which
  implies the lemma.
\end{proof}

\begin{theorem}
  \label{thm:cswalls}
  Let~$P$ be a 0/1-polytope such that each pair $s,t\in\verts{P}$,
  $s\not= t$, of distinct vertices~$s$ and~$t$ is a pair of antipodal
  cube vertices in a uniquely edge-cube spanned wall of~$P$. Then
  $\expans{\graph{P}}\geq 1$ holds.
\end{theorem}

\begin{proof}
  In each affine edge-cube spanning a uniquely edge-cube spanned wall
  of~$P$ we construct a flow as described in
  Observation~\ref{obs:cswalls}.  Let~$\phi$ be the sum of all these
  flows. Since each pair $s,t\in\verts{P}$, $s\not= t$, is a pair 
  of antipodal cube vertices in a uniquely edge-cube spanned wall
  of~$P$, the flow~$\phi$ has the properties required in
  Subsection~\ref{subsec:expflow}.  Lemma~\ref{lem:numwalls} ensures
  that each arc $(u,v)$ in the network~$\net{\graph{P}}$ is a cube-arc
  in at most~$\frac{n}{2}$ uniquely edge-cube spanned walls, if~$n$ is
  the number of vertices of~$P$. Thus we
  have~$\phi_{\max}\leq\frac{n}{2}$, and by~(\ref{eq:expphimax}) we
  obtain the claim of the theorem.
\end{proof}

Theorem~\ref{thm:cswalls} in particular yields a unified proof for the
following results which appeared in~\cite{Mih92} (where only a proof
for the statement concerning the perfect matching polytope is given).

\begin{corollary}
  \label{cor:cswalls}
  The graphs of the stable set polytope, the matching polytope, and
  the perfect matching polytope associated with an arbitrary graph have
  edge expansion at least one.
\end{corollary}

\begin{proof}
  Let~$G=(V,E)$ be a graph and let~$P$ be its stable set
  polytope. For two vertices~$s$ and~$t$ of~$P$
  let~$A_s,A_t\subseteq V$ be the corresponding stable sets in~$G$,
  and denote by $A^{(1)},\dots,A^{(k)}\subseteq V$ the node sets of
  the connected components of the subgraph of~$G$ induced by the
  symmetric difference of~$A_s$ and~$A_t$. Define
  $A^{(i)}_s:=A^{(i)}\cap A_s$ and $A^{(i)}_t:=A^{(i)}\cap A_t$.
  For each $\epsilon\in\{s,t\}^k$ the set
  $$
  S_{\epsilon}:=(A_s\cap A_t)\cup A^{(1)}_{\epsilon_1}\cup\cdots\cup A^{(k)}_{\epsilon_k}
  $$
  is stable in~$G$. By Chv\'{a}tal's result~\cite{Chv75} two
  vertices of~$P$ are adjacent if (and only if) the symmetric
  difference of the corresponding stable sets induces a connected
  subgraph of~$G$. Thus, the set~$C$ of vertices of~$P$ corresponding
  to $\setdef{S_{\epsilon}}{\epsilon\in\{s,t\}^k}$ is an affine
  edge-cube in~$P$, spanning the wall~$W$ which is defined by the
  equations $x_v=0$, $v\in V\setminus (A_s\cup A_t)$ and $x_v=1$,
  $v\in A_s\cap A_t$. Clearly, $s$ and~$t$ are antipodal vertices
  of~$C$. Since all pairs of mirror images in~$W$ belong to~$C$ (these
  pairs correspond to bipartitions of the subgraph of~$G$ that is
  induced by $A^{(1)}\cup\cdots A^{(k)}$), $W$ is
  uniquely edge-cube spanned by~$C$. Thus, $\expans{P}\geq 1$ by
  Theorem~\ref{thm:cswalls}.

  Since the matching polytope of a graph~$G$ is the stable set polytope of 
  the line graph of~$G$ (having the edges of~$G$ as vertices, which
  are adjacent if and only if the corresponding edges of~$G$ have a
  common end node), the claim on matching polytopes follows.
  
  Perfect matching polytopes satisfy the requirements of
  Theorem~\ref{thm:cswalls} as well; they even have the property that
  each pair of vertices spans a wall which \emph{is} an affine
  edge-cube.
\end{proof}

\section{Some Remarks}
\label{sec:remarks}

The results presented in this paper support the conjecture that graphs
of 0/1-polytopes inherently have good expansion properties and
therefore may in principle be good candidates for defining
neighborhood structures in the context of random walks. In fact, we
have proved for some classes of 0/1-polytopes, including simple
0/1-polytopes, stable set polytopes, and all 0/1-polytopes up to
dimension five, that their graphs have edge expansion at least one.

A proof of the conjecture that the edge expansion of the graph of any
$d$-dimensional 0/1-polytopes is bounded by the reciprocal of a
polynomial in~$d$ would have important consequences, even if this was
proved only for uniform 0/1-polytopes. For instance, such a result
would imply that indeed the bases-exchange graphs of arbitrary
matroids have sufficiently large edge expansion in order to construct a
randomized approximate counting algorithm. In particular, this would
solve the open questions for randomized approximation algorithms
for counting connected spanning subgraphs of a graph, forests of a
prescribed size in a graph, or maximal independent subsets in a given
set of vectors over $\GF{2}$ (see~\cite{JS97}).

Therefore, one might hope that, while the concept of the graph of a
0/1-polytope has not proven to be very useful in the context of
combinatorial optimization, it might have a successful revival in the
context of random generation and counting of certain combinatorial
objects.

\backmatter

\bibliographystyle{plain}
\bibliography{01expansion}

\end{document}